\documentclass[bezier,amstex, oneside,reqno]{amsart}
\usepackage{amsmath, amssymb, amsthm, verbatim}
\usepackage[dvips]{graphicx}
\usepackage{epsfig}
\usepackage{color}
\begingroup\makeatletter\ifx\SetFigFont\undefined%
\gdef\SetFigFont#1#2#3#4#5{%
  \reset@font\fontsize{#1}{#2pt}%
  \fontfamily{#3}\fontseries{#4}\fontshape{#5}%
  \selectfont}%
\fi\endgroup%

\newtheorem{theorem}{Theorem}
\newtheorem{lemma}{Lemma}
\newtheorem{prop}{Proposition}
\newtheorem{claim}{Claim}
\newtheorem{step}{Step}
\newtheorem*{reglemma}{Regularity Lemma}
\newtheorem*{livshitz}{Livshitz Theorem}
\newtheorem*{BI}{Theorem}

\theoremstyle{remark}
\newtheorem*{remark}{Remark}

\begin{document}
\author{A. Gogolev, M. Guysinsky}
\title[Smooth conjugacy of Anosov systems]{$C^1$-differentiable conjugacy of Anosov diffeomorphisms on three dimensional torus}
\dedicatory { Dedicated to Yakov Pesin on occasion of his 60th birthday.}
\begin{abstract}
We consider two $C^2$ Anosov diffeomorphisms in a $C^1$ neighborhood of a linear hyperbolic automorphism of three dimensional torus with real spectrum. We prove that they are $C^{1+\nu}$ conjugate if and only if the differentials of the return maps at corresponding periodic points have the same eigenvalues.
\end{abstract}
 \maketitle

\section{Introduction}
Consider an Anosov diffeomorphism $f$ of a compact smooth manifold. Structural stability asserts that if a diffeomorphism $g$ is $C^1$ close to $f$ then $f$ and $g$ are topologically conjugate. The conjugacy $h$ is unique in the homotopy class of identity.
$$h\circ f=g\circ h$$
It is known that $h$ is H\"older continuous.

There are simple obstructions for $h$ being smooth. Namely, let $x$ be a periodic point of $f$, $f^p(x)=x$. Then
$g^p(h(x))=h(x)$ and if $h$ were differentiable then
$$
Df^p(x)=\left(Dh(x)\right)^{-1}Dg^p(h(x))Dh(x)
$$
 i.e. $Df^p(x)$ and $Dg^p(h(x))$ are conjugate. We see that
every periodic point carries a modulus of $C^1$-differentiable conjugacy.

Suppose that for every periodic point $x$, $f^p(x)=x$, differentials of return maps $Df^p(x)$ and $Dg^p(h(x))$
are conjugate then we say that periodic data (p.~d.) of $f$ and $g$ coincide.

Suppose that p.~d. coincide, is $h$ differentiable?

A positive answer for Anosov diffeomorphisms of $\mathbb T^2$ was given in~\cite{LMM},~\cite{L}. De la Llave~\cite{L} observed that the answer is negative for Anosov diffeomorphisms of $\mathbb T^d$, $d\ge4$. He constructed two diffeomorphisms with the same p.~d. which are only H\"older conjugate.

We provide positive answer to the previous question in dimension three under an extra assumption.

The authors would like to thank A.Katok for suggesting us the
problem, numerous discussions and constant encouragement.

\section{Formulation of the main result}

Let $f$ be an Anosov diffeomorphism of $\mathbb T^d$. It is known~\cite{Manning} that $f$ is topologically
conjugate to a linear torus automorphism $L$. It is also known that Anosov diffeomorphisms of $\mathbb T^3$ are
the only Anosov diffeomorphisms on three dimensional manifolds~\cite{Franks},~\cite{Newhouse}.

Let $L$ be a hyperbolic automorphism of $\mathbb T^3$ with real eigenvalues. It is easy to show that absolute values of these eigenvalues are distinct. For the sake of notation we also assume that the eigenvalues are positive. This is not restrictive.

We will always assume that the Anosov diffeomorphisms that we are dealing with are at least $C^2$.

\begin{theorem} Given $L$ as above there exists a $C^1$-neighborhood $\mathcal U$ of $L$ such that any $f$ and $g$ in $\mathcal U$ having the same p.~d. are $C^{1+\nu}$ conjugate, $\nu>0$.
\end{theorem}

\begin{remark} The constant $\nu$ depends on the size of $\mathcal U$ and provided sufficient smoothness of $f$ and $g$ can be made as close as desired to $\log \lambda_3 / \log\lambda_2$ (see the definition in the next section) by shrinking the size of $\mathcal U$.
\end{remark}

\begin{remark} We don't know how to bootstrap regularity of $h$ to the regularity $f$ and $g$ like it was done in dimension two.
\end{remark}

A result about integrability of central distribution~\cite{BI} allows to show a stronger statement.
\begin{theorem}
Let $f$ and $g$ be Anosov diffeomorphisms of $\mathbb T^3$ and
$$h\circ f=g\circ h,$$
where $h$ is a homeomorphism homotopic to identity. Suppose that p.~d. coincide.

Also assume that $f$ and $g$ can be viewed as partially hyperbolic diffeomorphisms: there is an $f$-invariant splitting $T\mathbb T^3=E_f^s \oplus E_f^{wu} \oplus E_f^{su}$ and constants $c>0$, $0<\alpha<1<\tilde\beta<\beta<\gamma$ such that for $n>0$
\begin{multline}
\label{phd}
\|D(f^n)(x)(v)\|\le c\alpha^n\|v\|,\;\;\; v\in E_f^s(x),\\
\shoveleft{\frac1c\tilde\beta^n\|v\|\le\|D(f^n)(x)(v)\|\le c\beta^n\|v\|,\;\;\; v\in E_f^{wu}(x),}\\
\shoveleft{\frac1c\gamma^n\|v\|\le\|D(f^n)(x)(v)\|,\;\;\; v\in E_f^{su}(x)\hfill}
\end{multline}
Analogous conditions with possibly different set of constants hold for a $g$-invariant splitting $T\mathbb T^3=E_g^s \oplus E_g^{wu} \oplus E_g^{su}$.

Then the conjugacy $h$ is $C^{1+\nu}$, $\nu>0$.
\end{theorem}

\begin{remark} Here and further in the paper we assume that the unstable distribution has dimension two. Obviously one can formulate the counterpart of Theorem 2 in the case when stable distribution has dimension two.
\end{remark}

\section{Scheme of the proof}

Here we outline the proof of Theorem 1.

Let $\lambda_1$, $\lambda_2$ and $\lambda_3$ be the eigenvalues of the linear automorphism $L$, $0<\lambda_1<1<\lambda_2<\lambda_3$. We choose $\mathcal U$ in such a way that every $f\in\mathcal U$ is partially hyperbolic, satifying~(\ref{phd}) with constants $\alpha, \tilde\beta, \beta, \gamma$ independent on the choice of $f$, $0<\lambda_1<\alpha<1<\tilde\beta<\lambda_2<\beta<\gamma<\lambda_3$ and
\begin{equation}
\label{angles}
\measuredangle(E_L^\sigma, E_f^\sigma)<k<\frac\pi2, \; \sigma=s, wu, su.
\end{equation}

First we concentrate on a single diffeomorphism $f$ in $\mathcal U$. It is well known that distributions $E_f^s$, $E_f^u=E_f^{wu} \oplus E_f^{su}$ and $E_f^{su}$ integrate uniquely to stable, unstable and strong unstable foliations $W_f^s$, $W_f^u$ and $W_f^{su}$ respectively. We denote by $W_f^\sigma(x)$ the leaf of $W_f^\sigma$ passing through $x$, $\sigma=s, u, su$ and later $wu$. By $W_f^{\sigma}(x, R)$ we denote the local leaf of size $R$, i.~e., a ball of radius $R$ inside of $W_f^\sigma(x)$ centered at $x$, $\sigma=s, u, wu, su$. Let $h_f$ be conjugacy between $f$ and $L$, $h_f\circ f=L\circ h_f$. Stable and unstable foliations can be characterized topologically, e.g.
$$W_f^s(x)=\{y: d(f^n(x),f^n(y))\to 0\;\mbox{as}\;n\to+\infty\}.$$
As a consequence we have that $h_f(W_f^s)=W_L^s$ and $h_f(W_f^u)=W_L^u$. In other words $h_f$ maps leaves of foliations for $f$ into leaves of corresponding foliations for $L$.

We prove two simple lemmas.
\begin{lemma}
Let $f$ be in $\mathcal U$. Then the distribution $E_f^{wu}$ integrates uniquely to the foliation $W_f^{wu}$.
\end{lemma}
\begin{lemma} Define $h_f$ as above. Then
$h_f(W_f^{wu})=W_L^{wu}$.
\end{lemma}

Now let $f$ and $g$ be as in Theorem 1. For each of them we have the system of one dimensional invariant foliations. We know that $h(W_f^s)=W_g^s$. Also from Lemma 2 we have $h(W_f^{wu})=W_g^{wu}$ since $h=h_g^{-1}\circ h_f$. Consider restrictions of $h$ to the leaves of $W_f^s$ and $W_f^{wu}$. These restrictions are one dimensional maps. We show that they are smooth.

\begin{lemma} The conjugacy $h$ is $C^{1+\nu}$ along $W_f^s$.
\end{lemma}
 Which means that $h$ is differentiable along the stable foliation and the derivative is a H\"older continuous function on $\mathbb T^3$ with exponent $\nu$

\begin{remark}The general strategy of the proof of Theorem 1 is similar to de la Llave's strategy for Anosov diffeomorphisms of $\mathbb T^2$~\cite{L}. One proves smoothness of $h$ along one dimensional stable and unstable foliations. In particular proof of Lemma 3 can be carried out in the same way as in dimension two. The hard part is showing smoothness of $h$ along two dimensional unstable foliation.
\end{remark}
We would like to show the same for the foliation $W_f^{wu}$ but we split the proof into two steps.

\begin{lemma} The conjugacy $h$ is uniformly Lipschitz along $W_f^{wu}$.
\end{lemma}

\begin{lemma} The conjugacy $h$ is $C^{1+\nu}$ along $W_f^{wu}$.
\end{lemma}

After that we deal with the remaining foliation.

\begin{lemma} $h(W_f^{su})=W_g^{su}$.
\end{lemma}
\begin{remark}
We would like to remark that Lemma 6 requires only the coincidence of p.~d. in the weak unstable direction. It is not true in general that strong unstable foliations match.
\end{remark}

\begin{lemma} The conjugacy $h$ is $C^{1+\nu}$ along $W_f^{su}$.
\end{lemma}

\begin{remark} Proofs of smoothness along the foliations $W_f^{s}$ and $W_f^{su}$  are similar and use the coincidence of periodic data in corresponding directions. Showing smoothness along the weak unstable foliation is more subtle.\end{remark}

Now smoothness of $h$ is a simple consequence of a regularity result.

\begin{reglemma}\cite{J}  Let $M_j$ be a manifold and $W_j^s$, $W_j^u$ be continuous transverse foliations with uniformly smooth leaves, $j=1, 2$. Suppose that $h:M_1\to M_2$ is a homeomorphism that maps $W_1^s$ into $W_2^s$ and $W_1^u$ into $W_2^u$. Moreover assume that the restrictions of $h$ to the leaves of these foliations are uniformly $C^{r+\nu}, r\in \mathbb N,\; 0<\nu<1$, then $h$ is $C^{r+\nu}$.
\end{reglemma}

First we apply the lemma on every unstable leaf of $W_f^u$ for the pair of foliations $W_f^{wu}, W_f^{su}$. After we know that $h$ is $C^{1+\nu}$ along $W_f^u$ we finish by applying the lemma to stable and unstable foliations.

The structure of the next chapter is the following. We prove Lemmas 1 and 2 in Section 4.1. Section 4.2 is
devoted to the proof of Lemma 4. Sections 4.3 and 4.4 are the heart of our argument and contain proofs of Lemmas
5 and 6 correspondingly.


\section{Proof of Theorems 1 and 2}

First we prove Theorem 1.

\subsection{Weak unstable foliation}

In the proofs of Lemmas 1 and 2 we work with lifts of maps, distributions and foliations to $\mathbb R^3$. We
use the same notation for the lifts as for the objects themselves.

Denote by $d(\cdot, \cdot)$ the usual distance in $\mathbb R^3$ and let $d_f^\sigma(\cdot, \cdot)$ be the distance in the leaves of $W_f^\sigma$ which is defined only for pairs of points lying in the same leaf of $W_f^\sigma$, $\sigma =s, u, su, wu$.

\begin{proof}[Proof of Lemma 1]
Let us reason by contradiction. If $E_f^{wu}$ is not uniquely integrable then it must branch and we can find points $a, b, c \in \mathbb R^3$ such that
\begin{enumerate}
\item $a, b \in W_f^u(c)$,
\item there are smooth curves $\tau_{ca}, \tau_{cb}:[0, 1]\to W_f^u(c)$ such that $\tau_{ca}(0)=\tau_{cb}(0)=c$, $\tau_{ca}(1)=a$, $\tau_{cb}(1)=b$, and $\{\dot\tau_{ca}, \dot\tau_{cb}\}\subset E_f^{wu}$,
\item $a \in W_f^{su}(b)$.
\end{enumerate}
Then for $n\ge 1$
\begin{equation}
\label{dist1}
d(f^n(a), f^n(b))\le d(f^n(a), f^n(c))+d(f^n(c), f^n(b))\le c_1\beta^n,
\end{equation}
on the other hand

\begin{equation}
\label{dist2}
d_f^{su}(f^n(a), f^n(b))\ge c_2\gamma^n.
\end{equation}

For every $x\in\mathbb R^3$ consider a cone $Cone(x)=\{v\in T_x\mathbb R^3:  \measuredangle(v,E_L^{su}(x))\le
k\}$. The assumption~(\ref{angles}) tells us that $E_f^{su}(x)\subset Cone(x)$. Hence a leaf of $W_f^{su}$ can
be considered as a graph of a Lipschitz function over $E_L^{su}$. The Lipschitz constant depends only on $k$. It
follows that $W_f^{su}$ is quasi-isometric:

\begin{equation}
\label{QI}
\exists c_3>0\;\mbox{such that for}\; x\in W_f^{su}(y)\;\; d_f^{su}(x,y)\le c_3 d(x,y).
\end{equation}
Inequalities~(\ref{dist1}),~(\ref{dist2}) and~(\ref{QI}) sum up to a contradiction.
\end{proof}

\begin{proof}[Proof of Lemma 2]
Suppose that there are two points $a$ and $b$, $a \in W_f^{wu}(b)$ such that $h_f(a)\notin W_L^{wu}(h_f(b))$ then we have

\begin{equation}
\label{d1} d(f^n(a), f^n(b))\le c_1\beta^n,
\end{equation}
and since $h_f(a)$ and $h_f(b)$ lie in the same unstable leaf but not in the same weak unstable leaf we get
\begin{equation}
\label{d2}
d(h_f(f^n(a)),h_f(f^n(b)))=d(L^n(h_f(a)),L^n(h_f(b)))\ge c_2\gamma^n.
\end{equation}
Finally since
\begin{equation}
\label{periodicity}
h_f(x+\bar m)=h_f(x)+\bar m,\;\;\bar m\in \mathbb Z^3
\end{equation}
we have that $d(h_f(x),h_f(y))\le c(\varepsilon) d(x,y)$ for any $x$ and $y$ such that $d(x,y)\ge\varepsilon$. Hence

\begin{equation}
\label{d3}
d(h_f(f^n(a)), h_f(f^n(b)))\le c_3d(f^n(a), f^n(b)),
\end{equation}
where $c_3$ depends on $d(a, b)$.
Inequalities~(\ref{d1}),~(\ref{d2}) and (\ref{d3}) sum up to a contradiction.
\end{proof}


\subsection{Affine structure on the weak unstable foliation}

Let $f$ be in $\mathcal U$. For any $x$ and $y$, $y\in W_f^{wu}(x)$ define the function
$$\rho_f(x, y)=\prod_{n\ge 1}\frac{D_f^{wu}(f^{-n}(y))}{D_f^{wu}(f^{-n}(x))}$$
where $D_f^{wu}(z)=\|D(f)\big|_{E_f^{wu}}(z)\|$.
The following properties are easy to prove:
\begin{itemize}
\item[(P1)] $\rho_f(x,\cdot)$ is well defined and H\"older continuous.
\item[(P2)] $\forall x, y \in W_f^{wu}(z)\;\; \rho_f(x, y)\rho_f(y, z)=\rho_f(x, z)$.
\item[(P3)] $\rho_f(f(x),f(y))=\frac{D_f^{wu}(y)}{D_f^{wu}(x)}\rho_f(x, y)$.
\item[(P4)] The function $\rho(\cdot, \cdot)$ is the only continuous function satisfying $\rho_f(x,x)=1$ and Property 3.
\item[(P5)] $\forall K>0\;\; \exists C>0$ such that $C>\rho_f(x, y)>\frac1C$ whenever $d^{wu}(x, y)<K$.
\end{itemize}

The goal is to show that $h$ is differentiable along $W_f^{wu}$ ($wu$-differentiable) and

\begin{equation}
\label{hrho}
\rho_g(h(x),h(y))=\frac{D_h^{wu}(y)}{D_h^{wu}(x)}\rho_f(x,y),
\end{equation}

\begin{proof}[Proof of Lemma 4] Fix an arbitrary point $p$. Let $h_p: W_f^{wu}(p)\to W_f^{wu}(h(p))$ be the restriction of $h$ to $W_f^{wu}(p)$. We would like to show that $h_p$ is Lipschitz with a constant that does not depend on $p$. Let $m$ be the induced volume on $W_f^{wu}(p)$. Consider the function $\tilde d_f$
$$\tilde d_f(x,y)=\int\limits_x^y\frac{1}{\rho_f(x,z)}dm(z),\;\; x,y\in W_f^{wu}(p),$$
we integrate along the leaf with respect to the measure $m$.

Function $\tilde d_f$ has the following properties which are simple corollories of the properties of $\rho_f$ and the definition of $\tilde d_f$.

\begin{itemize}
\item[(D1)] $\tilde d_f(x,y)=d_f^{wu}(x,y)+o(d_f^{wu}(x,y))$,
\item[(D2)] $\tilde d_f(f(x),f(y))=D_f^{wu}(x)\tilde d_f(x,y)$,
\item[(D3)] $\forall K>0 \;\exists C>0$ such that
\begin{equation}
\frac1C\tilde d_f(x,y)\le d_f^{wu}(x,y)\le C\tilde d_f(x,y)
\end{equation}
whenever $d_f^{wu}(x,y)<K$.
\item[(D4)] The function $\tilde d_f$ is continuous. To state this property precisely we consider lift of $\tilde d_f$. We speak about lifts of points and leaves.
\begin{multline*}
\forall\varepsilon>0\; \exists\delta>0\;\;\;\mbox{such that}\;\;\; \forall x,y\in \mathbb R^3,
y\in W_f^{wu}(x)\;\;\;\\
\mbox{and}\;\;\; \forall z,q\in\mathbb R^3, q\in W_f^{wu}(z), z\in B(x,\delta), q\in B(y,\delta)\\
\;\;\mbox{we have}\;\;\;\; |\tilde d_f(x,y)-\tilde d_f(z,q)|<\varepsilon.
\end{multline*}
\end{itemize}
We will also need $\tilde d_g$ which is defined analogously on the leaves of $W_g^{wu}$ and has analogous properties.

The lift of the conjugacy $h$ satisfies the equation~(\ref{periodicity}) which implies the following
$$\exists C>0: \forall x,y \;\;\; d(h(x),h(y))\le C d(x,y)\;\;\mbox{if}\;\;d(x,y)\ge 1.$$
Also we know that weak unstable foliation is quasi-isometric which gives us the same for the distance in weak unstable foliations
\begin{equation}
\label{farlipschitz}
\exists C>0: \forall x,y \;\;\; d_g^{wu}(h(x),h(y))\le C d_f^{wu}(x,y)\;\;\mbox{if}\;\;d_f^{wu}(x,y)\ge 1.
\end{equation}
This tells us that $h_p$ is Lipschitz for points that are far enough. So we need to estimate $d_g^{wu}(h(x),h(y))$ for $x$ and $y$ close. Note that~(D3) allows us to use $\tilde d_g$ and $\tilde d_f$ in these estimates instead of $d_g^{wu}$ and $d_f^{wu}$.

Recall the following well-known result.
\begin{livshitz}
If $f:M\to M$ is a transitive Anosov diffeomorphism and $\varphi_1, \varphi_2:M\to\mathbb R$ are H\"older continuous functions such that
$$\prod_{i=1}^p\varphi_1(f^i(x))=\prod_{i=1}^p\varphi_2(f^i(x))\;\;\mbox{whenever}\;f^p(x)=x$$
then there is a function $P:M\to\mathbb R$, unique up to a multiplicative constant, such that
$$\frac{\varphi_1}{\varphi_2}=\frac{P\circ f}{P}.$$
Moreover $P$ is H\"older continuous.
\end{livshitz}

Apply Livshitz Theorem for $\varphi_1=D_f^{wu}(\cdot)$ and $\varphi_2=D_g^{wu}(h(\cdot))$. The condition of the
Livshitz Theorem is satisfied because of the assumption on p.~d. We have
\begin{equation}
 \forall n>0 \;\;\;\;\;
\prod_{i=0}^{n-1}\frac{D_g^{wu}\left(h(f^i(x))\right)}{D_f^{wu}(f^i(x))}=\frac{P(x)}{P(f^n(x))}.
\end{equation}

Choose points $x$ and $y$ close on the leaf $W_f^{wu}(p)$. Choose the smallest $N$ such that $d_f^{wu}(f^N(x),f^N(y))\ge 1$. Then
\begin{multline*}
\frac{\tilde d_g(h(x),h(y))}{\tilde d_f(x,y)}=
\prod_{i=0}^{N-1}\frac{D_g^{wu}\left(g^i(h(x))\right)}{D_f^{wu}(f^i(x))}\cdot
\frac {\tilde d_g\left(g^N(h(x)),g^N(h(y))\right)} {\tilde d_f(f^N(x),f^N(y))}\\
=\frac{P(x)}{P(f^N(x))}\cdot\frac {\tilde d_g\left(g^N(h(x)),g^N(h(y))\right)} {\tilde d_f(f^N(x),f^N(y))}\le
\frac{P(x)}{P(f^N(x))}\cdot constant.
\end{multline*}
Here we used~(\ref{farlipschitz}) and (D3) for $\tilde d_f$ and $\tilde d_g$.
Function $P$ is bounded away from zero and infinity so we get that $h$ is uniformly Lipschitz along the weak unstable foliation.
\end{proof}

\subsection{Transitive point argument and construction of a measure absolutely continuous with respect to weak unstable foliation}

We divide the proof of Lemma 5 into several steps. The conjugacy $h$ is Lipschitz along $W_f^{wu}$ and hence $wu$-differentiable at almost every point with respect to Lebesgue measure on the leaves of $W_f^{wu}$. It is obvious that $wu$-differentiability of $h$ at $x$ implies $wu$-differentiability of $h$ at any point from the orbit $\{f^i(x), i\in \mathbb Z\}$. Moreover:
\begin{step} Suppose that $h$ is wu-differentiable at $x$ and $\overline{\{f^i(x), i\ge 0\}}=\mathbb T^3$ then $h$ is $C^{1+\nu}$ along $W_f^{wu}$ and~(\ref{hrho}).
\end{step}

The problem now is to show existence of such a transitive point $x$. We know that almost every point is
transitive with respect to a given ergodic measure with full support. On the other hand $h$ is
$wu$-differentiable at almost every point with respect to Lebesgue measure on the leaves. Unfortunately it can
happen that for natural ergodic "physical measures" these two "full measure" sets do not intersect. In other
words weak unstable foliation is not absolutely continuous with respect to a "physical measure".

Let us explain this phenomenon in more detail. Consider a volume preserving $C^1$ small perturbation $\tilde L$ of $L$, $H\circ L=\tilde L\circ H$. The Lyapunov exponents of $\tilde L$ are defined on a full volume set of regular points $\mathcal R$ and are given by the formula
$$\chi^\sigma=\int_{\mathbb T^3}\log D_{\tilde L}^\sigma d \mbox{vol},\;\;\; \sigma=s, wu, su.$$
The perturbation $\tilde L$ can be chosen in such a way that $\chi^{wu}>\log\lambda_2$ (see~\cite{BB}, Proposition 0.3). It is easy to show that the weak unstable foliation of $\tilde L$ is not absolutely continuous. Namely, let $\Delta$ be a segment of a weak unstable leaf of $L$. Then by Lemma 2 $H(\Delta)$ is a piece of a weak unstable leaf of $\tilde L$. We show that Lebesgue measure of $\mathcal R\cap H(\Delta)$ is equal to zero. For any $n\ge 0$ $H(L^n(\Delta))=\tilde L^n(H(\Delta))$ and~(\ref{angles}) guarantees that $\tilde L^n(H(\Delta))$ can be viewed as a graph of a Lipschitz function over a leaf of the weak unstable foliation of $L$. Hence
$$\mbox{length}(\tilde L^n(H(\Delta))\le c_1\cdot \mbox{length}(L^n(\Delta))=\lambda_2^n\cdot \mbox{length}(\Delta),\;\; n\ge 0.$$
Suppose that $\mbox{Leb}(\mathcal R\cap H(\Delta))>0$ then
$$\mbox{length}(\tilde L^n(H(\Delta))\ge c_2e^{n(\chi^{wu}-\varepsilon)},\;\;\varepsilon=\frac12(\chi^{wu}-\log\lambda_2)$$
which contradicts the previous inequality.

This observation answers a question of Hirayama and Pesin~\cite{HP} about existence of non-absolutely continuous foliations with non-compact leaves.

To overcome this problem we do
\begin{step}
Construction of a measure $\mu$ absolutely continuous with respect to $W_f^{wu}$.
\end{step}
This construction follows the lines of Pesin-Sinai~\cite{PS} construction of $u$-Gibbs measures. In our setup
the construction is simpler so for the sake of completeness we present it here. Measure $\mu$ has full support.
Thus ergodicity of $\mu$ would imply that almost every point is transitive and hence by Step 1 $h$ would be
$wu$-differentiable. We do not know how to show ergodicity of $\mu$. Instead we do
\begin{step}
Set of transitive points is a full measure $\mu$ set.
\end{step}
Steps 2 and 3 guarantee existence of a transitive point needed in Step 1.

\begin{proof}[Proof of Lemma 5]
~\\
{\bfseries Step 1.} Let us pick a point $y\in \mathbb T^3$ and show that $h$ is $wu$-differentiable at $y$ and moreover
\begin{equation}
\label{wudiff}
D_h^{wu}(y)=\frac{P(y)}{P(x)}D_h^{wu}(x)
\end{equation}
where $P$ is the same as in the proof of Lemma 4.

Choose $y'\in W_f^{wu}(y)$. Property (D1) of $\tilde d_f$, $\tilde d_g$ ensures that  it is enough to show that
\begin{equation}
\label{wutilde}
\frac{\tilde d_g(h(y),h(y'))}{\tilde d_f(y,y')}=\frac{P(y)}{P(x)}D_h^{wu}(x).
\end{equation}

Fix an $\varepsilon >0$ small compared to $\tilde d_f(y,y')$. Choose a small open ball $B$ centered at $y$ and
define

\begin{multline*}
B'=\{z' : \;\;\exists z\in B\;\;\mbox{such that}\;\; \tilde d_f(z,z')=\tilde d_f(y,y')\\
 \mbox{and} \;\; (z,z') \;\;\mbox{has the same orientation as}\;\; (y,y')\}.
\end{multline*}
The condition about orientation ensures that $B'$ has only one connected component. The set $B'$ is a small neighborhood of $y'$ because of the continuity of $\tilde d_f$~(D4). The size of $B$ must be chosen in such a way that
\begin{enumerate}
\item $|P(z)-P(y)|<\varepsilon$ if $z\in B$,
\item $|\tilde d_g(h(z),h(z'))-\tilde d_g(h(y),h(y'))|<\varepsilon$ where $z$ and $z'$ are the same as in definition of $B'$.
\end{enumerate}

Since $x$ is transitive there is an arbitrarily large $N$ such that $f^N(x)\in B$. Choose $z$ on $W_f^{wu}(x)$ such that $\tilde d_f(f^N(x),f^N(z))=\tilde d_f(y,y')$ so that $f^N(z)\in B'$ by the definition. We choose $N$ big enough so that
$$\left|\frac{\tilde d_g(h(x),h(z))}{\tilde d_f(x,z)}-D_h^{wu}(x)\right|<\varepsilon.$$

\begin{figure}[htbp]
\begin{center}
\scalebox{0.7}{
\begin{picture}(0,0)%
\epsfig{file=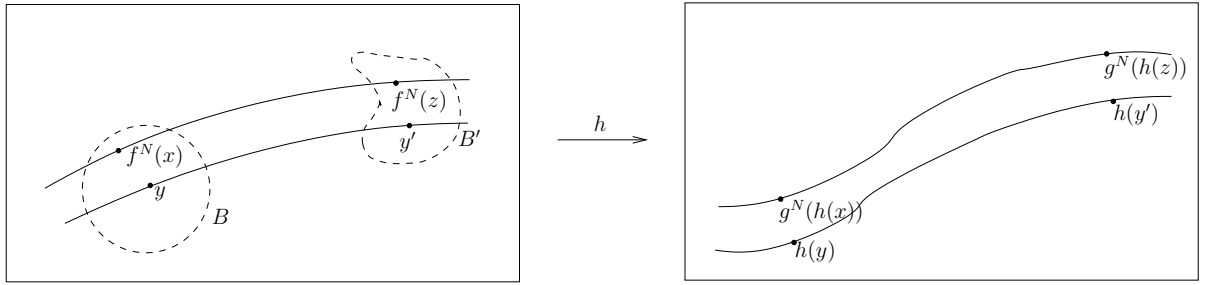}%
\end{picture}%
\setlength{\unitlength}{3947sp}%
\begingroup\makeatletter\ifx\SetFigFont\undefined%
\gdef\SetFigFont#1#2#3#4#5{%
  \reset@font\fontsize{#1}{#2pt}%
  \fontfamily{#3}\fontseries{#4}\fontshape{#5}%
  \selectfont}%
\fi\endgroup%
\begin{picture}(10719,2529)(857,-2371)
\put(2206,-1579){\makebox(0,0)[lb]{\smash{\SetFigFont{12}{14.4}{\familydefault}{\mddefault}{\updefault}{\color[rgb]{0,0,0}$y$}%
}}}
\put(1951,-1272){\makebox(0,0)[lb]{\smash{\SetFigFont{12}{14.4}{\familydefault}{\mddefault}{\updefault}{\color[rgb]{0,0,0}$f^N(x)$}%
}}}
\put(4404,-1151){\makebox(0,0)[lb]{\smash{\SetFigFont{12}{14.4}{\familydefault}{\mddefault}{\updefault}{\color[rgb]{0,0,0}$y'$}%
}}}
\put(4329,-783){\makebox(0,0)[lb]{\smash{\SetFigFont{12}{14.4}{\familydefault}{\mddefault}{\updefault}{\color[rgb]{0,0,0}$f^N(z)$}%
}}}
\put(7952,-2126){\makebox(0,0)[lb]{\smash{\SetFigFont{12}{14.4}{\familydefault}{\mddefault}{\updefault}{\color[rgb]{0,0,0}$h(y)$}%
}}}
\put(7809,-1774){\makebox(0,0)[lb]{\smash{\SetFigFont{12}{14.4}{\familydefault}{\mddefault}{\updefault}{\color[rgb]{0,0,0}$g^N(h(x))$}%
}}}
\put(10809,-896){\makebox(0,0)[lb]{\smash{\SetFigFont{12}{14.4}{\familydefault}{\mddefault}{\updefault}{\color[rgb]{0,0,0}$h(y')$}%
}}}
\put(10727,-484){\makebox(0,0)[lb]{\smash{\SetFigFont{12}{14.4}{\familydefault}{\mddefault}{\updefault}{\color[rgb]{0,0,0}$g^N(h(z))$}%
}}}
\put(6147,-994){\makebox(0,0)[lb]{\smash{\SetFigFont{12}{14.4}{\familydefault}{\mddefault}{\updefault}{\color[rgb]{0,0,0}$h$}%
}}}
\put(2719,-1834){\makebox(0,0)[lb]{\smash{\SetFigFont{12}{14.4}{\familydefault}{\mddefault}{\updefault}{\color[rgb]{0,0,0}$B$}%
}}}
\put(4932,-1136){\makebox(0,0)[lb]{\smash{\SetFigFont{12}{14.4}{\familydefault}{\mddefault}{\updefault}{\color[rgb]{0,0,0}$B'$}%
}}}
\end{picture}
}
\end{center}
\caption{Differentiability of $h$ at the point $y$.}
\label{fig1}
\end{figure}
Now we are ready to do the estimates
\begin{multline*}
\tilde d_g(h(y),h(y'))=\varepsilon_1+\tilde d_g\left(h(f^N(x)),h(f^N(z))\right)\\
=\varepsilon_1+\tilde d_g\left(g^N(h(x)),g^N(h(z))\right)\\
=\varepsilon_1+\prod\nolimits_{i=0}^{N-1}D_g^{wu}\left(g^i(h(x))\right)\cdot\tilde d_g\left(h(x),h(z)\right)\\
=\varepsilon_1+\prod\nolimits_{i=0}^{N-1}D_g^{wu}\left(h(f^i(x))\right) (D_h^{wu}(x)+\varepsilon_2)\tilde d_f(x,z)\\
=\varepsilon_1+\frac{\prod\nolimits_{i=0}^{N-1}D_g^{wu}\left(h(f^i(x))\right)}{\prod\nolimits_{i=0}^{N-1}D_f^{wu}(f^i(x))} (D_h^{wu}(x)+\varepsilon_2)\tilde d_f(f^N(x),f^N(z))\\
=\varepsilon_1+\frac{P(f^N(x))}{P(x)} (D_h^{wu}(x)+\varepsilon_2)\tilde d_f(y, y')\\
=\varepsilon_1+\frac{P(y)+\varepsilon_3}{P(x)} (D_h^{wu}(x)+\varepsilon_2)\tilde d_f(y, y')
\end{multline*}
with $\max(|\varepsilon_1|, |\varepsilon_2|, |\varepsilon_3|)\le \varepsilon$. Now letting $\varepsilon$ go to $0$ we get~(\ref{wutilde}).

To show~(\ref{hrho}) define

$$\tilde \rho_g(h(x),h(y))=\frac{D_h^{wu}(y)}{D_h^{wu}(x)}\rho_f(x,y).$$
Then

\begin{multline*}
\tilde \rho_g\left(g(h(x)),g(h(y))\right)=\tilde \rho_g\left(h(f(x)),h(f(y))\right)\\
=\frac{D_h^{wu}(f(y))}{D_h^{wu}(f(x))}\rho_f(f(x),f(y))=\frac{D_h^{wu}(f(y))}{D_h^{wu}(f(x))}\cdot\frac{D_f^{wu}(y)}{D_f^{wu}(x)}\rho_f(x,y)\\
=\frac{D_{h\circ f}^{wu}(y)}{D_{h\circ f}^{wu}(x)} \rho_f(x,y) =\frac{D_{g\circ h}^{wu}(y)}{D_{g\circ
h}^{wu}(x)} \rho_f(x,y) =\frac{D_g^{wu}(h(y))}{D_g^{wu}(h(x))}\cdot\frac{D_h^{wu}(y)}{D_h^{wu}(x)}\rho_f(x,y)=\\
\frac{D_g^{wu}(h(y))}{D_g^{wu}(h(x))}\tilde\rho_g(h(x),h(y)).
\end{multline*}
This by $(P4)$ implies that $\tilde \rho_g=\rho_g$ which is equivalent to~(\ref{hrho}).

{\bfseries Step 2.}  Let $x_0$ be a fixed point for $f$ and let $V_0$ be an open bounded neighborhood of $x_0$ in $W_f^{wu}(x_0)$. Consider a probability measure $\eta^0$ supported on $V_0$ with density proportional to $\rho_f(x_0,\cdot)$. For $n>0$ define
$$V_n=f^n(V_0),\; \eta^n=(f^n)_*\eta^0$$
so that $\eta^n$ is supported on $V_n$ and has density proportional to $\rho_f(x_0,\cdot)$ by $(P3)$.

Let $\mu^n=\frac1n\sum_{i=0}^{n-1}\eta^i$. By the Krylov-Bogoljubov theorem $\{\mu^n; n\ge 0\}$ is weakly compact and any of its limits is $f$-invariant. Let $\mu$ be a one of those limits along a subsequence $\{n_k; k\ge 1\}$. We would like to prove that $\mu$ has absolutely continuous conditional measures on the pieces of weak unstable foliation.

Let us be more precise. Consider a small open set $X\subset\mathbb T^3$ which can be decomposed in the following way
$$X=\bigcup_{y\in Y}W_f^{wu}(y,R_y).$$
Here $Y$ is a two dimensional transversal. To simplify the notation let $W(y)=W_f^{wu}(y,R_y)$. Denote by $\mu_T$ the transverse measure on $Y$: for $Y'\subset Y$ $\mu_T(Y')=\mu(\cup_{y\in Y'}W(y))$. Similary define $\eta_T^n$ and $\mu_T^n$. Obviously $\mu_T^{n_k}\to\mu_T$ weakly as $k\to\infty$. We show that for $\mu_T$ almost every $y$, $y\in Y$ the conditional measure $\mu_y$ on the local leaf $W(y)$ is absolutely continuous with respect to Lebesgue measure $m_y$ on $W(y)$.

The conditional measures are characterized by the following property

\begin{equation}
\label{rohlin}
\forall F\in C(X)\;\; \int_X F d\mu=\int_Y d\mu_T(y)\int_{W(y)}F(y,z)d\mu_y(z).
\end{equation}

First we look at conditional measures of $\eta^n$. We fix $X$ and $Y$ as above and we assume that the end points of $V_n$ lie outside of $X$. Let $\{a_1, a_2,\ldots a_m\}=Y\cap V_n$. Then the formulas for the transverse measure and conditional measures are obvious:
\begin{multline}
\label{cond}
\eta_T^n=\sum_{i=1}^m\left(\int_{W(a_i)}\rho_f(x_0,z)dm_{a_i}(z)\right)\delta(a_i),\\
\shoveleft{d\eta_y^n(z)}=\left(\int_{W(y)}\rho_f(y,z)dm_y(z)\right)^{-1}\rho_f(y,z)dm_y(z).\hfill
\end{multline}
Notice that $\eta_y^n$ actually do not depend on $n$.

The goal now is to show that $d\mu_y=\left(\int_{W(y)}\rho_f(y,z)dm_y(z)\right)^{-1}\rho_f(y,\cdot)$ for almost every $y$. It could happen that the end points of $V_n$ lie inside of $X$. Support $S_n$ of $\eta_T^n$ consists of finitely many points. Some of these points correspond to the end points of $V_n$. Denote the set of these points by $B_n$, $|B_n|\in\{0,1,2\}$. Let $A_n=S_n\backslash B_n$ then there is a natural decomposition of the transverse measure $\eta_T^n$


\begin{figure}[htbp]
\begin{center}
\scalebox{0.7}{
\begin{picture}(0,0)%
\epsfig{file=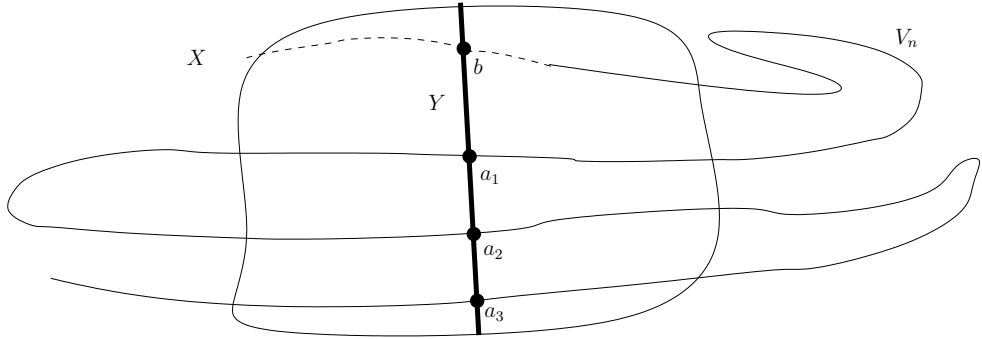}%
\end{picture}%
\setlength{\unitlength}{3947sp}%
\begingroup\makeatletter\ifx\SetFigFont\undefined%
\gdef\SetFigFont#1#2#3#4#5{%
  \reset@font\fontsize{#1}{#2pt}%
  \fontfamily{#3}\fontseries{#4}\fontshape{#5}%
  \selectfont}%
\fi\endgroup%
\begin{picture}(8744,3073)(1211,-2658)
\put(5412,-267){\makebox(0,0)[lb]{\smash{\SetFigFont{12}{14.4}{\familydefault}{\mddefault}{\updefault}{\color[rgb]{0,0,0}$b$}%
}}}
\put(5472,-1219){\makebox(0,0)[lb]{\smash{\SetFigFont{12}{14.4}{\familydefault}{\mddefault}{\updefault}{\color[rgb]{0,0,0}$a_1$}%
}}}
\put(5502,-1894){\makebox(0,0)[lb]{\smash{\SetFigFont{12}{14.4}{\familydefault}{\mddefault}{\updefault}{\color[rgb]{0,0,0}$a_2$}%
}}}
\put(5510,-2449){\makebox(0,0)[lb]{\smash{\SetFigFont{12}{14.4}{\familydefault}{\mddefault}{\updefault}{\color[rgb]{0,0,0}$a_3$}%
}}}
\put(2840,-184){\makebox(0,0)[lb]{\smash{\SetFigFont{12}{14.4}{\familydefault}{\mddefault}{\updefault}{\color[rgb]{0,0,0}$X$}%
}}}
\put(5009,-589){\makebox(0,0)[lb]{\smash{\SetFigFont{12}{14.4}{\familydefault}{\mddefault}{\updefault}{\color[rgb]{0,0,0}$Y$}%
}}}
\put(9194,  4){\makebox(0,0)[lb]{\smash{\SetFigFont{12}{14.4}{\familydefault}{\mddefault}{\updefault}{\color[rgb]{0,0,0}$V_n$}%
}}}
\end{picture}

}
\end{center}
\caption{Decomposition of the transverse measure $\nu_T^n$.}
\label{fig2}
\end{figure}
\begin{multline*}
\eta_T^n=\sum_{a\in A_n}\left(\int_{W(a)}\rho(x_0,z)dm_a(z)\right)\delta(a)\\
+\sum_{b\in B_n}\left(\int_{W(b)\cap V_n}\rho(x_0,z)dm_b(z)\right)\delta(b)=\eta_{(T,A)}^n+\eta_{(T,B)}^n.
\end{multline*}
The conditional measures $\eta_y^n$ for $y\notin B_n$ are given by formula~(\ref{cond}).
Since $W_f^{wu}$ is uniformly expanding it is clear that $\eta_T^n(B_n)\to 0$ as $n\to\infty$. Hence
\begin{equation}
\label{endpoints}
\frac1n\left(\sum_{i=0}^{n-1}\eta_T^i(B_i)\right)\to 0,\;\;\;n\to\infty
\end{equation}
and

\begin{equation}
\label{muA}
\mu_T=\lim_{k\to\infty}\mu_T^{n_k}=\lim_{k\to\infty}\frac1{n_k}\sum_{i=0}^{n_k-1}\eta_{(T,A)}^i.
\end{equation}

Consider a continuous function $F$ on $X$.

\begin{multline*}
\int_XFd\mu=\lim_{k\to\infty}\int_XFd\mu^{n_k}=\lim_{k\to\infty}\frac1{n_k}\sum_{i=0}^{n_k-1}\int_XFd\eta^{n_k}\\
=\lim_{k\to\infty}\frac1{n_k}\sum_{i=0}^{n_k-1}\int_Yd\eta_T^{n_k}(y)\int_{W(y)}F(y,z)d\eta_k^{n_k}(z)\\
=\lim_{k\to\infty}\frac1{n_k}\sum_{i=0}^{n_k-1}\int_Yd\eta_{(T,A)}^{n_k}(y)\int_{W(y)}F(y,z)d\eta_y^{n_k}(z)\\
+\lim_{k\to\infty}\frac1{n_k}\sum_{i=0}^{n_k-1}\int_Yd\eta_{(T,B)}^{n_k}(y)\int_{W(y)}F(y,z)d\eta_y^{n_k}(z).
\end{multline*}
The function $F$ is bounded so it follows from~(\ref{endpoints}) that the last limit is zero. So we get

\begin{multline*}
\int_XFd\mu=\lim_{k\to\infty}\frac1{n_k}\sum_{i=0}^{n_k-1}\int_Yd\eta_{(T,A)}^{n_k}(y)\int_{W(y)}F(y,z)d\eta_y^{n_k}(z)\\
=\lim_{k\to\infty}\frac1{n_k}\sum_{i=0}^{n_k-1}\int_Yd\eta_{(T,A)}^{n_k}(y)\left(\int_{W(y)}\rho_f(y,z)dm_y(z)\right)^{-1}\int_{W(y)}F(y,z)\rho_f(y,z)dm_y(z).
\end{multline*}
Now notice that the function that we integrate with respect to $\eta_{(T,A)}^{n_k}$ is continuous and does not depend on $n_k$. Hence using~(\ref{muA}) we get

\begin{multline*}
\int_XFd\mu=\lim_{k\to\infty}\int_Yd\mu_T^{n_k}(y)\left(\int_{W(y)}\rho_f(y,z)dm_y(z)\right)^{-1}\int_{W(y)}F(y,z)\rho_f(y,z)dm_y(z)\\
=\int_Yd\mu_T(y)\left(\int_{W(y)}\rho_f(y,z)dm_y(z)\right)^{-1}\int_{W(y)}F(y,z)\rho_f(y,z)dm_y(z)
\end{multline*}
and by~(\ref{rohlin}) we see that up to normalization the density of the conditional measure on $W(y)$ is equal to $\rho_f(y,\cdot)$ for $\mu_T$ a.~e. $y$.

The leaf $W_f^{wu}(x_0)$ is dense in $\mathbb T^3$ since $W_f^{wu}(x_0)=h_f^{-1}\left(W_L^{wu}(h_f^{-1}(x_0))\right)$ and $W_L^{wu}(h_f^{-1}(x_0))$ is a dense irrational line in $\mathbb T^3$.
Hence the support $\mu$ is the whole torus.

{\bfseries Step 3.}
To prove that $\mu$ a.~e. point is transitive we fix a ball in $\mathbb T^3$ and show that a.~e. point visits the ball infinitely many times. Then to conclude transitivity we only need to cover $\mathbb T^3$ by a countable collection of balls such that every point is contained in an arbitrarily small ball.

So let us fix a ball $B'$ and a slightly smaller ball $B$, $B\subset B'$. Let $\psi$ be a non-negative continuous function supported on $B'$  and equal to $1$ on $B$. By Birkhoff ergodic theorem

\begin{equation}
\label{ergodictheorem}
E(\psi|\mathcal I)=\lim_{n\to\infty}\frac1n\sum_{i=0}^{n-1}\psi\circ f^i
\end{equation}
where $\mathcal I$ is $\sigma$-algebra of $f$-invariant sets.

Let $A=\{x: E(\psi|\mathcal I)(x)=0\}$. Then $\mu(A\cap B)=0$ since $\int_A\psi d\mu=\int_A E(\psi|\mathcal I)d\mu=0$. Hence

\begin{equation*}
\label{positiveonB}
E(\psi|\mathcal I)(x)>0\;\; \mbox{for} \;\;\mu \;\mbox{a.e.}\;\; x\in B.
\end{equation*}

Since $h_f(W_f^{wu})=W_L^{wu}$ it is possible to find $R>0$ such that $\cup_{b\in B}W(b, R)=\mathbb T^3$.

\begin{remark}
This observation also implies that $\mu$ has full support.
\end{remark}
Applying the standard Hopf argument we get that for $\mu$ a.~e. $x$ the function $E(\psi|\mathcal I)$ is constant on $W(x, R)$. Now absolute continuity of $W_f^{wu}$ together with above observations shows that $E(\psi|\mathcal I)>0$ for $\mu$ a.~e. $x$ which means according to~(\ref{ergodictheorem}) that a.~e. $x$ visits $B'$ infinitely many times.
\end{proof}

\subsection{Strong unstable foliations match}

Let us point out once again that in the proof of Lemma 6 we only use $wu$-differentiability of $h$ which as we showed is equivalent to coincidence of p.~d. in the weak unstable direction.

\begin{proof}[Proof of Lemma 6]

We will be working on two dimensional leaves of $W_f^u$. We know that each of these leaves is subfoliated by $W_f^{wu}$ as well as by $W_f^{su}$. The goal is to prove that $h(W_f^{su})=W_g^{su}$ so we consider the foliation $U=h^{-1}(W_g^{su})$. As for usual foliations $U(x)$ stands for the leaf of $U$ passing through $x$ and $U(x,R)$ stands for the local leaf of size $R$. Obviously $U$ subfoliate $W_f^u$. A priori the leaves of $U$ are just H\"older continuous curves. Since weak unstable foliations match we see that a leaf $U(x)$ intersects each $W_f^{wu}(y)$, $y\in W_f^u(x)$ exactly once.

Let us prove several auxiliary claims.

\begin{claim} Consider a point $a\in\mathbb T^3$. Suppose that there is a point $b\ne a$, $b\in W_f^{su}(a)\cap U(a)$. Let $c\in W_f^{wu}(a)$ and $d=W_f^{wu}(b)\cap W_f^{su}(c)$, $e=W_f^{wu}(b)\cap U(c)$. Then $d=e$.
\end{claim}

\begin{figure}[h!]
\begin{center}
\scalebox{0.8}{
\begin{picture}(0,0)%
\epsfig{file=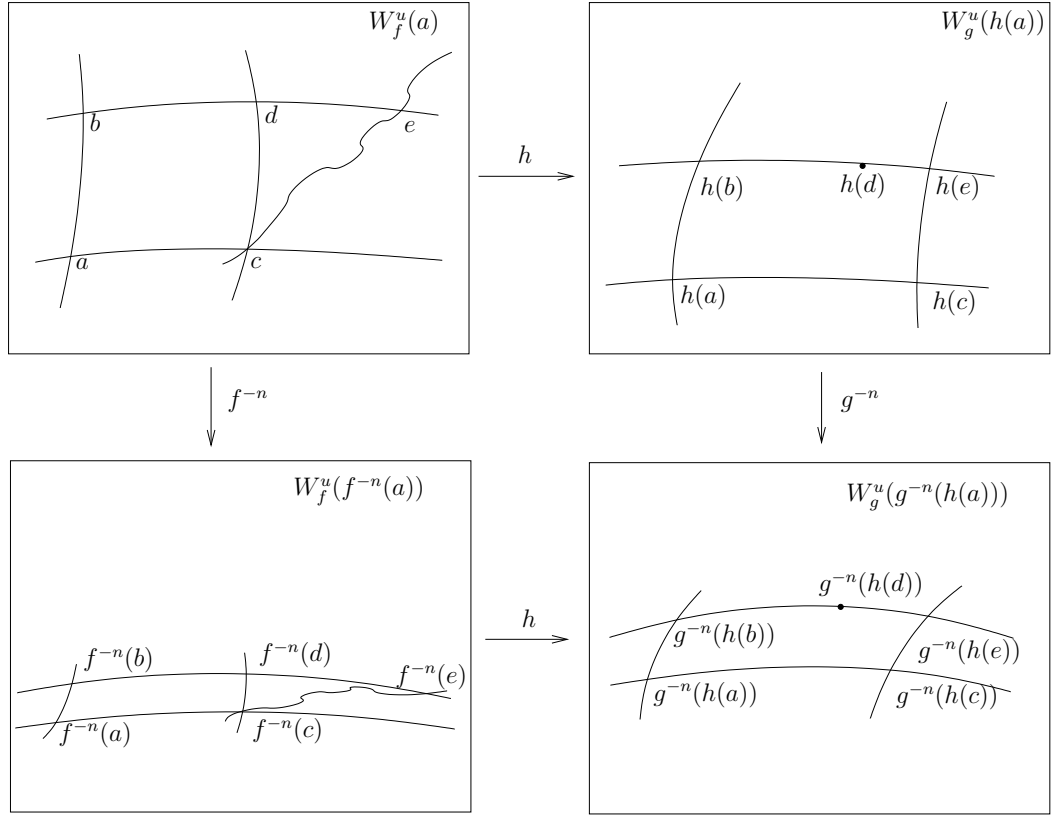}%
\end{picture}%
\setlength{\unitlength}{3947sp}%
\begingroup\makeatletter\ifx\SetFigFont\undefined%
\gdef\SetFigFont#1#2#3#4#5{%
  \reset@font\fontsize{#1}{#2pt}%
  \fontfamily{#3}\fontseries{#4}\fontshape{#5}%
  \selectfont}%
\fi\endgroup%
\begin{picture}(8199,6399)(1808,-6076)
\put(2353,-1781){\makebox(0,0)[lb]{\smash{\SetFigFont{12}{14.4}{\familydefault}{\mddefault}{\updefault}{\color[rgb]{0,0,0}$a$}%
}}}
\put(2465,-701){\makebox(0,0)[lb]{\smash{\SetFigFont{12}{14.4}{\familydefault}{\mddefault}{\updefault}{\color[rgb]{0,0,0}$b$}%
}}}
\put(4933,-709){\makebox(0,0)[lb]{\smash{\SetFigFont{12}{14.4}{\familydefault}{\mddefault}{\updefault}{\color[rgb]{0,0,0}$e$}%
}}}
\put(4655,101){\makebox(0,0)[lb]{\smash{\SetFigFont{12}{14.4}{\familydefault}{\mddefault}{\updefault}{\color[rgb]{0,0,0}$W_f^u(a)$}%
}}}
\put(9148,101){\makebox(0,0)[lb]{\smash{\SetFigFont{12}{14.4}{\familydefault}{\mddefault}{\updefault}{\color[rgb]{0,0,0}$W_g^u(h(a))$}%
}}}
\put(7258,-1196){\makebox(0,0)[lb]{\smash{\SetFigFont{12}{14.4}{\familydefault}{\mddefault}{\updefault}{\color[rgb]{0,0,0}$h(b)$}%
}}}
\put(7100,-2044){\makebox(0,0)[lb]{\smash{\SetFigFont{12}{14.4}{\familydefault}{\mddefault}{\updefault}{\color[rgb]{0,0,0}$h(a)$}%
}}}
\put(2218,-5487){\makebox(0,0)[lb]{\smash{\SetFigFont{12}{14.4}{\familydefault}{\mddefault}{\updefault}{\color[rgb]{0,0,0}$f^{-n}(a)$}%
}}}
\put(3740,-5457){\makebox(0,0)[lb]{\smash{\SetFigFont{12}{14.4}{\familydefault}{\mddefault}{\updefault}{\color[rgb]{0,0,0}$f^{-n}(c)$}%
}}}
\put(2413,-4902){\makebox(0,0)[lb]{\smash{\SetFigFont{12}{14.4}{\familydefault}{\mddefault}{\updefault}{\color[rgb]{0,0,0}$f^{-n}(b)$}%
}}}
\put(3792,-4864){\makebox(0,0)[lb]{\smash{\SetFigFont{12}{14.4}{\familydefault}{\mddefault}{\updefault}{\color[rgb]{0,0,0}$f^{-n}(d)$}%
}}}
\put(4063,-3559){\makebox(0,0)[lb]{\smash{\SetFigFont{12}{14.4}{\familydefault}{\mddefault}{\updefault}{\color[rgb]{0,0,0}$W_f^u(f^{-n}(a))$}%
}}}
\put(4865,-5030){\makebox(0,0)[lb]{\smash{\SetFigFont{12}{14.4}{\familydefault}{\mddefault}{\updefault}{\color[rgb]{0,0,0}$f^{-n}(e)$}%
}}}
\put(8405,-3604){\makebox(0,0)[lb]{\smash{\SetFigFont{12}{14.4}{\familydefault}{\mddefault}{\updefault}{\color[rgb]{0,0,0}$W_g^u(g^{-n}(h(a)))$}%
}}}
\put(6890,-5172){\makebox(0,0)[lb]{\smash{\SetFigFont{12}{14.4}{\familydefault}{\mddefault}{\updefault}{\color[rgb]{0,0,0}$g^{-n}(h(a))$}%
}}}
\put(8795,-5194){\makebox(0,0)[lb]{\smash{\SetFigFont{12}{14.4}{\familydefault}{\mddefault}{\updefault}{\color[rgb]{0,0,0}$g^{-n}(h(c))$}%
}}}
\put(7048,-4714){\makebox(0,0)[lb]{\smash{\SetFigFont{12}{14.4}{\familydefault}{\mddefault}{\updefault}{\color[rgb]{0,0,0}$g^{-n}(h(b))$}%
}}}
\put(8188,-4309){\makebox(0,0)[lb]{\smash{\SetFigFont{12}{14.4}{\familydefault}{\mddefault}{\updefault}{\color[rgb]{0,0,0}$g^{-n}(h(d))$}%
}}}
\put(8967,-4841){\makebox(0,0)[lb]{\smash{\SetFigFont{12}{14.4}{\familydefault}{\mddefault}{\updefault}{\color[rgb]{0,0,0}$g^{-n}(h(e))$}%
}}}
\put(5825,-949){\makebox(0,0)[lb]{\smash{\SetFigFont{12}{14.4}{\familydefault}{\mddefault}{\updefault}{\color[rgb]{0,0,0}$h$}%
}}}
\put(5855,-4594){\makebox(0,0)[lb]{\smash{\SetFigFont{12}{14.4}{\familydefault}{\mddefault}{\updefault}{\color[rgb]{0,0,0}$h$}%
}}}
\put(3545,-2869){\makebox(0,0)[lb]{\smash{\SetFigFont{12}{14.4}{\familydefault}{\mddefault}{\updefault}{\color[rgb]{0,0,0}$f^{-n}$}%
}}}
\put(8360,-2869){\makebox(0,0)[lb]{\smash{\SetFigFont{12}{14.4}{\familydefault}{\mddefault}{\updefault}{\color[rgb]{0,0,0}$g^{-n}$}%
}}}
\put(9103,-1189){\makebox(0,0)[lb]{\smash{\SetFigFont{12}{14.4}{\familydefault}{\mddefault}{\updefault}{\color[rgb]{0,0,0}$h(e)$}%
}}}
\put(9080,-2096){\makebox(0,0)[lb]{\smash{\SetFigFont{12}{14.4}{\familydefault}{\mddefault}{\updefault}{\color[rgb]{0,0,0}$h(c)$}%
}}}
\put(8352,-1174){\makebox(0,0)[lb]{\smash{\SetFigFont{12}{14.4}{\familydefault}{\mddefault}{\updefault}{\color[rgb]{0,0,0}$h(d)$}%
}}}
\put(3725,-1774){\makebox(0,0)[lb]{\smash{\SetFigFont{12}{14.4}{\familydefault}{\mddefault}{\updefault}{\color[rgb]{0,0,0}$c$}%
}}}
\put(3838,-634){\makebox(0,0)[lb]{\smash{\SetFigFont{12}{14.4}{\familydefault}{\mddefault}{\updefault}{\color[rgb]{0,0,0}$d$}%
}}}
\end{picture}
}
\end{center}
\caption{Illustration to the proof of the Claim 1. Notice that the actual size of the bottom pictures should be much smaller.}
\label{fig3}
\end{figure}
Assume that $d\ne e$. For the sake of concreteness we also assume that $d$ lies between $b$ and $e$. We look at
configurations $\{a,b,c,d,e\}\in W_f^u(a)$, $\{h(a),h(b)$, $h(c)$,$h(d)$,$h(e)\}$ $\in W_g^u(h(a))$ and study
their evolution under $f^{-n}, n>0$ and $g^{-n}, n>0$ respectively. Since under the action of $f^{-1}$ strong
unstable leaves contract exponentially faster then weak unstable leaves we get that

\begin{equation}
\label{acbd}
\forall \varepsilon>0 \;\; \exists n_0: \;\forall n>n_0\;\;\;
\left| \frac{d_f^{wu}(f^{-n}(a),f^{-n}(c))}{d_f^{wu}(f^{-n}(b),f^{-n}(d))}-1 \right|<\varepsilon.
\end{equation}
Analogously

\begin{equation}
\label{Hacbd}
\forall \varepsilon>0 \;\; \exists n_1: \;\forall n>n_1\;\;\;
\left| \frac{d_g^{wu}\left(g^{-n}(h(a)),g^{-n}(h(c))\right)}{d_g^{wu}\left(g^{-n}(h(b)),g^{-n}(h(e))\right)}-1 \right|<\varepsilon.
\end{equation}
The next statement is a direct corollary of~(D1) and~(D2). There exists a $\delta>0$ which depends on the initial configuration $\{a,b,c,d,e\}$ such that

\begin{equation}
\label{bebd}
\forall n>0 \;\;\;\;\; \frac{d_f^{wu}(f^{-n}(b),f^{-n}(e))}{d_f^{wu}(f^{-n}(b),f^{-n}(d))} >1+\delta.
\end{equation}
Combining~(\ref{acbd}) and~(\ref{bebd}) we get

\begin{equation}
\label{differentDistance}
\exists \delta'>0:\;\forall n>n_0 \;\;\;\;\; \frac{d_f^{wu}(f^{-n}(b),f^{-n}(e))}{d_f^{wu}(f^{-n}(a),f^{-n}(c))} >1+\delta'.
\end{equation}
On the other hand we know that $h$ is continuously $wu$-differentiable, hence

\begin{multline}
\label{sameDistance} \forall \varepsilon>0\;\; \exists n_2:\; \forall n>n_2\;\;\;\;\;
\left| \frac{d_g^{wu}\left(g^{-n}(h(a)), g^{-n}(h(c))\right)}{d_f^{wu}(f^{-n}(a), f^{-n}(c))}-D_h^{wu}(f^{-n}(a))\right|<\varepsilon\\
\mbox{and}\;\;  \left| \frac{d_g^{wu}\left(g^{-n}(h(b)), g^{-n}(h(e))\right)}{d_f^{wu}(f^{-n}(b),
f^{-n}(e))}-D_h^{wu}(f^{-n}(a))\right|<\varepsilon.
\end{multline}
It is easy to see that~(\ref{differentDistance}) contradicts~(\ref{sameDistance}) and~(\ref{Hacbd}) so we are
done.

\begin{claim}
Consider a weak unstable leaf $W_f^{wu}(a)$ and $b\in W_f^{su}(a)$, $b\ne a$. For any $y\in W_f^{wu}(a)$ let
$y'=W_f^{wu}(b)\cap W_f^{su}(y)$. Then $\exists c_1, c_2>0$ such that $\forall y\in W_f^{wu}(a)$
$\;\;\;c_1>d_f^{su}(y,y')>c_2$.
\end{claim}

Recall that $h_f(W_f^{wu})=W_L^{wu}$. The leaves $W_L^{wu}(h_f(a))$ and $W_L^{wu}(h_f(b))$ are parallel lines in
$W_L^u(h_f(a))$ that are fixed distance apart. Hence the estimate from below is a direct consequence of uniform
continuity of $h_f|_{W_f^u}(a)$ with respect to metrics $d_f^u$ and $d_L^u$.

Now we prove the estimate from above. We need to show that the strip between $W_f^{wu}(a)$ and $W_f^{wu}(b)$ cannot contain arbitrarily long pieces of strong unstable leaves. The reason for this is uniform transversality of weak unstable and strong unstable foliation.

For any positive number $\varepsilon$ we can choose a finite number of points $\{a_0=a, a_1, a_2,\ldots
a_m=b\}$ between $a$ and $b$ on $W_f^{su}(a)$ in such a way that $W_f^{wu}(a_i)$ is contained in
$\varepsilon$-neighborhood of $W_f^{wu}(a_{i-1})$ and vice versa, $i=1,\ldots m$. Again this is possible because
$W_L^{wu}(h_f(a_i))$, $i=0,\ldots m$ are parallel lines and $h_f^{-1}|_{W_L^u(h_f(a))}$ is uniformly continuous.

Let $s=\min_{x\in\mathbb T^3} \measuredangle(E_f^{wu}(x),E_f^{su}(x))$. Choose a small $\delta>0$ such that in
any ball $B$ of size $\delta$

$$\max_{x, y\in B}\max\{\measuredangle(E_f^{wu}(x),E_f^{wu}(y)), \measuredangle(E_f^{su}(x),E_f^{su}(y))\}<\frac {s}{10}.$$
In such a ball the direction of $E_f^{wu}$ is almost constant comparing to the angle between $E_f^{wu}$ and
$E_f^{su}$. Clearly it is possible to choose a small $\varepsilon=\varepsilon(s,\delta)$  and correspondingly
the points $\{a_0, a_1,\ldots a_m\}$ as above such that any strong unstable leave crosses the strip between
$W_f^{wu}(a_{i-1})$ and $W_f^{wu}(a_i)$ in a ball of size $\delta$, $i=1,\ldots m$.
 This gives us uniform estimates on the lengths of pieces of strong unstable leaves in the strips between  $W_f^{wu}(a_{i-1})$ and $W_f^{wu}(a_i)$, $i=\overline {1,m}$. The sum of these estimates gives us the desired uniform estimate from above.

\begin{claim}
Suppose $\exists a\in \mathbb T^3$ and $R>0$ such that $W_f^{su}(a, R)=U(a, R)$ then $W_f^{su}=U$.
\end{claim}

Consider a point $c\in W_f^{wu}(a)$ then applying Claim 1 to the points $b\in W_f^{su}(a, R)$ we get that
$\exists R_c>0$ such that $W_f^{su}(c, R_c)=U(c, R_c)$. Moreover by Claim 2 numbers $R_c$, $c\in W_f^{wu}(a)$
are uniformly bounded away from zero. Now the statement follows from denseness of $W_f^{wu}(a)$ in $\mathbb
T^3$.

We are ready to prove the lemma.

We say that $W_f^{su}(x)$ and $U(x)$ intersect transversally at $y$ if $y\in W_f^{su}(x)\cap U(x)$ and $\forall R>0\;$ the local leaf $U(y,R)$ lies on both sides of $W_f^{su}(y)$.

We consider two cases.

{\bfseries Case 1.} {\itshape At every periodic point $x_0$ the leaves $W_f^{su}(x_0)$ and $U(x_0)$ do not
intersect transversally at a point different from $x_0$.}

Notice that the property of having a transverse intersection is stable --- if $W_f^{su}(x)$ and $U(x)$ intersect transversally then there is a neighborhood $V$ of $x$ such that $\forall z\in V$ $W_f^{su}(z)$ and $U(z)$ intersect transversally. Periodic points are dense therefore absence of transverse intersections at periodic points leads to absence of transverse intersections at all points.

We assume that $W_f^{su}\ne U$. Then the above observation together with Claim 3 tell us that for any point $x$ the leaves $W_f^{su}(x)$ and $U(x)$ intersect only at $x$.
\begin{figure}[h!]
\begin{center}
\scalebox{0.8}{
\begin{picture}(0,0)%
\epsfig{file=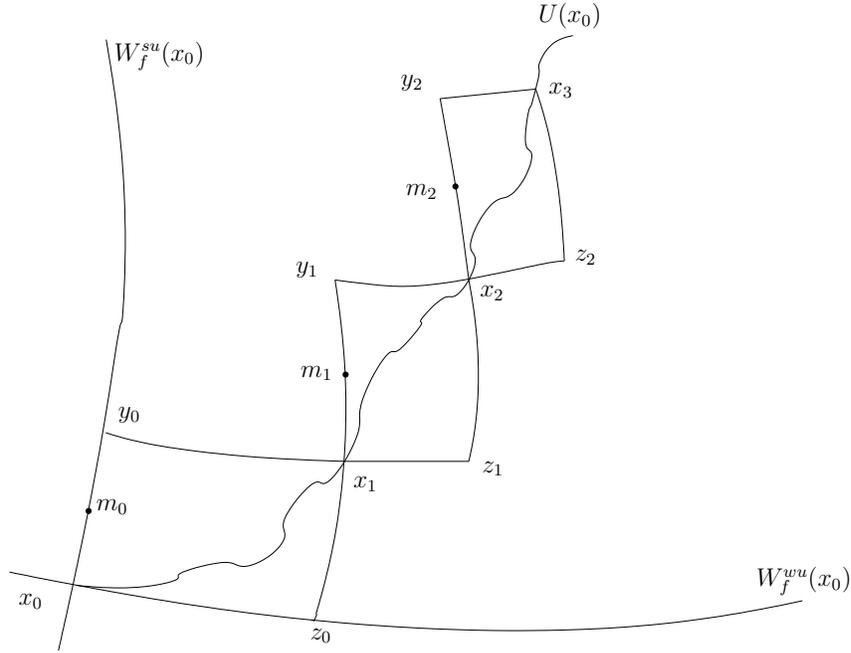}%
\end{picture}%
\setlength{\unitlength}{3947sp}%
\begingroup\makeatletter\ifx\SetFigFont\undefined%
\gdef\SetFigFont#1#2#3#4#5{%
  \reset@font\fontsize{#1}{#2pt}%
  \fontfamily{#3}\fontseries{#4}\fontshape{#5}%
  \selectfont}%
\fi\endgroup%
\begin{picture}(6249,5143)(2158,-6209)
\put(3001,-1561){\makebox(0,0)[lb]{\smash{\SetFigFont{12}{14.4}{\familydefault}{\mddefault}{\updefault}{\color[rgb]{0,0,0}$W_f^{su}(x_0)$}%
}}}
\put(2251,-5836){\makebox(0,0)[lb]{\smash{\SetFigFont{12}{14.4}{\familydefault}{\mddefault}{\updefault}{\color[rgb]{0,0,0}$x_0$}%
}}}
\put(2858,-5094){\makebox(0,0)[lb]{\smash{\SetFigFont{12}{14.4}{\familydefault}{\mddefault}{\updefault}{\color[rgb]{0,0,0}$m_0$}%
}}}
\put(3031,-4366){\makebox(0,0)[lb]{\smash{\SetFigFont{12}{14.4}{\familydefault}{\mddefault}{\updefault}{\color[rgb]{0,0,0}$y_0$}%
}}}
\put(4539,-6106){\makebox(0,0)[lb]{\smash{\SetFigFont{12}{14.4}{\familydefault}{\mddefault}{\updefault}{\color[rgb]{0,0,0}$z_0$}%
}}}
\put(4876,-4914){\makebox(0,0)[lb]{\smash{\SetFigFont{12}{14.4}{\familydefault}{\mddefault}{\updefault}{\color[rgb]{0,0,0}$x_1$}%
}}}
\put(5888,-4794){\makebox(0,0)[lb]{\smash{\SetFigFont{12}{14.4}{\familydefault}{\mddefault}{\updefault}{\color[rgb]{0,0,0}$z_1$}%
}}}
\put(5873,-3414){\makebox(0,0)[lb]{\smash{\SetFigFont{12}{14.4}{\familydefault}{\mddefault}{\updefault}{\color[rgb]{0,0,0}$x_2$}%
}}}
\put(6616,-3136){\makebox(0,0)[lb]{\smash{\SetFigFont{12}{14.4}{\familydefault}{\mddefault}{\updefault}{\color[rgb]{0,0,0}$z_2$}%
}}}
\put(6413,-1824){\makebox(0,0)[lb]{\smash{\SetFigFont{12}{14.4}{\familydefault}{\mddefault}{\updefault}{\color[rgb]{0,0,0}$x_3$}%
}}}
\put(6331,-1261){\makebox(0,0)[lb]{\smash{\SetFigFont{12}{14.4}{\familydefault}{\mddefault}{\updefault}{\color[rgb]{0,0,0}$U(x_0)$}%
}}}
\put(4426,-3226){\makebox(0,0)[lb]{\smash{\SetFigFont{12}{14.4}{\familydefault}{\mddefault}{\updefault}{\color[rgb]{0,0,0}$y_1$}%
}}}
\put(5250,-1756){\makebox(0,0)[lb]{\smash{\SetFigFont{12}{14.4}{\familydefault}{\mddefault}{\updefault}{\color[rgb]{0,0,0}$y_2$}%
}}}
\put(4464,-4044){\makebox(0,0)[lb]{\smash{\SetFigFont{12}{14.4}{\familydefault}{\mddefault}{\updefault}{\color[rgb]{0,0,0}$m_1$}%
}}}
\put(5289,-2618){\makebox(0,0)[lb]{\smash{\SetFigFont{12}{14.4}{\familydefault}{\mddefault}{\updefault}{\color[rgb]{0,0,0}$m_2$}%
}}}
\put(8041,-5686){\makebox(0,0)[lb]{\smash{\SetFigFont{12}{14.4}{\familydefault}{\mddefault}{\updefault}{\color[rgb]{0,0,0}$W_f^{wu}(x_0)$}%
}}}
\end{picture}

}
\end{center}
\caption{The ladder of rectangles.}
\label{fig4}
\end{figure}
Let $x_0$ to be a fixed point of $f$. For each $y\in U(x_0)$ the leaf $W_f^{su}(y)$ intersects $U(x_0)$ only at
$y$. Thus we are able to build a ladder of rectangles in $W_f^u(x_0)$ as shown on the Figure \ref{fig4}. The
sides of the rectangles are pieces of weak unstable and strong unstable leaves. The rectangles are subject to
condition
$$d_f^{su}(x_i, y_i)=1,\;i\ge 0.$$
This guarantees that after the choice of $y_0$ (there are two choices) the sequence of rectangles is defined
uniquely. Let $d_i=d_f^{wu}(y_i, x_{i+1}),\; i\ge 0$ and let $\{m_i;\;i\ge 0\}$ be midpoints on the sides of
rectangles as shown on the picture.

Suppose that $\inf_{i\ge 0}d_i>0$. Apply $f^{-n},\; n>0$ to the ladder of rectangles. The leaf $U(x_0)$ is invariant while the rectangles shrink and become flat as shown on Figure \ref{fig5}. Namely

$$\forall \varepsilon>0\;\;\exists n_0:\;\forall n>n_0\;\;\;\mbox{and}\;\;\forall i\ge 0\;\;
\frac{d_f^{su}(f^{-n}(x_i),f^{-n}(y_i))}{d_f^{wu}(f^{-n}(y_i),f^{-n}(x_{i+1}))}<\varepsilon.$$

\begin{figure}[htbp]
\begin{center}
\scalebox{0.8}{
\begin{picture}(0,0)%
\epsfig{file=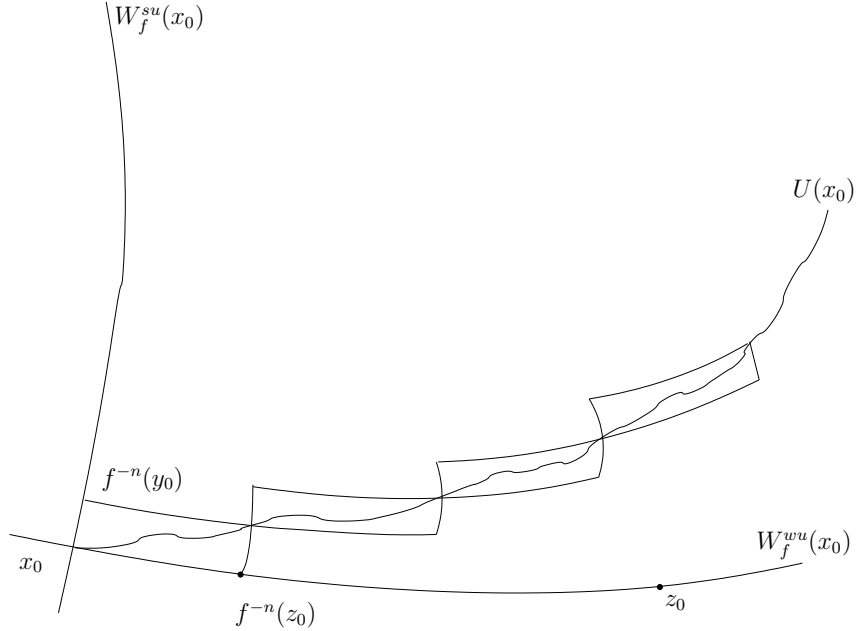}%
\end{picture}%
\setlength{\unitlength}{3947sp}%
\begingroup\makeatletter\ifx\SetFigFont\undefined%
\gdef\SetFigFont#1#2#3#4#5{%
  \reset@font\fontsize{#1}{#2pt}%
  \fontfamily{#3}\fontseries{#4}\fontshape{#5}%
  \selectfont}%
\fi\endgroup%
\begin{picture}(6450,4883)(2158,-6249)
\put(3001,-1561){\makebox(0,0)[lb]{\smash{\SetFigFont{12}{14.4}{\familydefault}{\mddefault}{\updefault}{\color[rgb]{0,0,0}$W_f^{su}(x_0)$}%
}}}
\put(2251,-5836){\makebox(0,0)[lb]{\smash{\SetFigFont{12}{14.4}{\familydefault}{\mddefault}{\updefault}{\color[rgb]{0,0,0}$x_0$}%
}}}
\put(8041,-5686){\makebox(0,0)[lb]{\smash{\SetFigFont{12}{14.4}{\familydefault}{\mddefault}{\updefault}{\color[rgb]{0,0,0}$W_f^{wu}(x_0)$}%
}}}
\put(7321,-6129){\makebox(0,0)[lb]{\smash{\SetFigFont{12}{14.4}{\familydefault}{\mddefault}{\updefault}{\color[rgb]{0,0,0}$z_0$}%
}}}
\put(2896,-5184){\makebox(0,0)[lb]{\smash{\SetFigFont{12}{14.4}{\familydefault}{\mddefault}{\updefault}{\color[rgb]{0,0,0}$f^{-n}(y_0)$}%
}}}
\put(3945,-6249){\makebox(0,0)[lb]{\smash{\SetFigFont{12}{14.4}{\familydefault}{\mddefault}{\updefault}{\color[rgb]{0,0,0}$f^{-n}(z_0)$}%
}}}
\put(8326,-2919){\makebox(0,0)[lb]{\smash{\SetFigFont{12}{14.4}{\familydefault}{\mddefault}{\updefault}{\color[rgb]{0,0,0}$U(x_0)$}%
}}}
\end{picture}
}
\end{center}
\caption{Ladder of rectangles after several iterations.}
\label{fig5}
\end{figure}
This means that in any fixed bounded neighborhood of $x_0$ the leaf
$U(x_0)$ is arbitrarily close to $W_f^{wu}(x_0)$. In particular we
have that $x_1$ is arbitrarily close to $z_0$ while we know that
they are some fixed distance apart. To make this argument completely
rigorous one needs to carry out an estimate on the distance between
$z_0$ and $x_1$ using regularity of holonomies along $W_f^{wu}$ and
$W_f^{su}$ inside of the leaf $W_f^u(x_0)$.
 We conclude that $\inf_{i\ge 0}d_i=0$.

Then choose a subsequence $\{m_{n_k}; k\ge 0\}$ such that corresponding rectangles have width going to zero as $k$ tend to infinity. Each of these rectangles contains a piece of $U(x_0)$ inside of it. Let $m$ be an accumulation point of $\{m_{n_k}; k\ge 0\}$ considered as a sequence of points in $\mathbb T^3$ rather than on $W_f^u(x_0)$. Since the width of the rectangles is shrinking and the foliations are continuous we get that $W_f^{su}(m,\frac12)=U(m,\frac12)$. Hence $W_f^{su}=U$ by Claim 3 and we move on to the second case.

{\bfseries Case 2.} {\itshape
There exist a periodic point $x_0$ and a point $y_0, y_0\ne x_0$ such that $W_f^{su}(x_0)$ and $U(x_0)$ intersect at $y_0$ transversally.}

Without loss of generality we can assume that $x_0$ is a fixed point. We chose a sequence $\{x_i\in W_f^{wu}(x_0); i\ge 1\}$ such that $x_i\to y_0$, $i\to\infty$. Here and afterwards we speak about convergence on the torus, not in the leaf $W_f^u(x_0)$. By Claim 1 we know that for any $i$ the leaves $W_f^{wu}(y_0)$, $W_f^{su}(x_i)$ and $U(x_i)$ intersect at one point $z_i$. Up to the choice of a subsequence we have that $z_i\to y_1$, $i\to\infty$, where $y_1$ is some point on $W_f^{su}(y_0)$. Since the foliation $U$ is continuous we have that $y_1\in U(y_0)=U(x_0)$ as well. The strong unstable foliation is orientable and the pairs $(x_0, y_0), (x_i, z_i), i\ge 1$ have the same orientation i.~e. $y_0$ lies between $x_0$ and $y_1$.
\begin{figure}[htbp]
\begin{center}
\scalebox{0.7}{
\begin{picture}(0,0)%
\epsfig{file=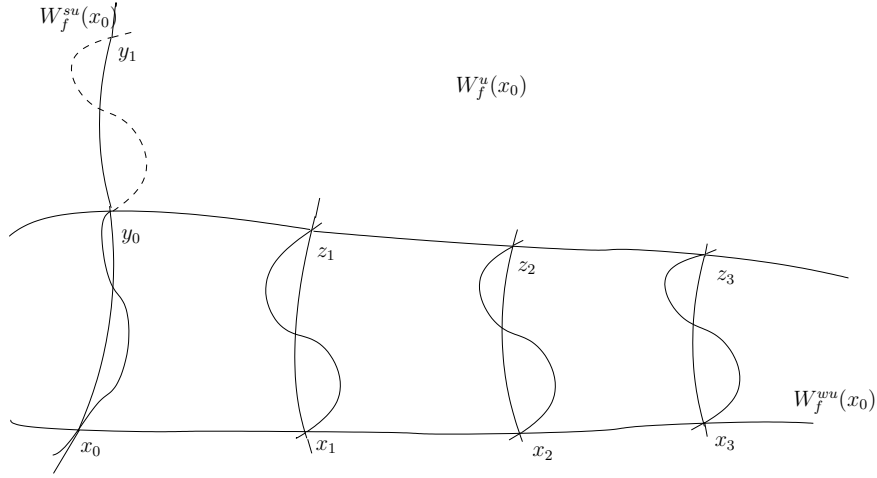}%
\end{picture}%
\setlength{\unitlength}{3947sp}%
\begingroup\makeatletter\ifx\SetFigFont\undefined%
\gdef\SetFigFont#1#2#3#4#5{%
  \reset@font\fontsize{#1}{#2pt}%
  \fontfamily{#3}\fontseries{#4}\fontshape{#5}%
  \selectfont}%
\fi\endgroup%
\begin{picture}(7554,4263)(1959,-4231)
\put(2616,-4017){\makebox(0,0)[lb]{\smash{\SetFigFont{12}{14.4}{\familydefault}{\mddefault}{\updefault}{\color[rgb]{0,0,0}$x_0$}%
}}}
\put(4731,-3994){\makebox(0,0)[lb]{\smash{\SetFigFont{12}{14.4}{\familydefault}{\mddefault}{\updefault}{\color[rgb]{0,0,0}$x_1$}%
}}}
\put(6674,-4054){\makebox(0,0)[lb]{\smash{\SetFigFont{12}{14.4}{\familydefault}{\mddefault}{\updefault}{\color[rgb]{0,0,0}$x_2$}%
}}}
\put(8309,-3964){\makebox(0,0)[lb]{\smash{\SetFigFont{12}{14.4}{\familydefault}{\mddefault}{\updefault}{\color[rgb]{0,0,0}$x_3$}%
}}}
\put(8324,-2471){\makebox(0,0)[lb]{\smash{\SetFigFont{12}{14.4}{\familydefault}{\mddefault}{\updefault}{\color[rgb]{0,0,0}$z_3$}%
}}}
\put(6539,-2389){\makebox(0,0)[lb]{\smash{\SetFigFont{12}{14.4}{\familydefault}{\mddefault}{\updefault}{\color[rgb]{0,0,0}$z_2$}%
}}}
\put(4739,-2261){\makebox(0,0)[lb]{\smash{\SetFigFont{12}{14.4}{\familydefault}{\mddefault}{\updefault}{\color[rgb]{0,0,0}$z_1$}%
}}}
\put(2984,-2104){\makebox(0,0)[lb]{\smash{\SetFigFont{12}{14.4}{\familydefault}{\mddefault}{\updefault}{\color[rgb]{0,0,0}$y_0$}%
}}}
\put(2954,-469){\makebox(0,0)[lb]{\smash{\SetFigFont{12}{14.4}{\familydefault}{\mddefault}{\updefault}{\color[rgb]{0,0,0}$y_1$}%
}}}
\put(9015,-3598){\makebox(0,0)[lb]{\smash{\SetFigFont{12}{14.4}{\familydefault}{\mddefault}{\updefault}{\color[rgb]{0,0,0}$W_f^{wu}(x_0)$}%
}}}
\put(5985,-778){\makebox(0,0)[lb]{\smash{\SetFigFont{12}{14.4}{\familydefault}{\mddefault}{\updefault}{\color[rgb]{0,0,0}$W_f^u(x_0)$}%
}}}
\put(2250,-163){\makebox(0,0)[lb]{\smash{\SetFigFont{12}{14.4}{\familydefault}{\mddefault}{\updefault}{\color[rgb]{0,0,0}$W_f^{su}(x_0)$}%
}}}
\end{picture}

}
\end{center}
\caption{Curves $U(x_i)$ that pass through $x_i$ and $z_i$ are the preimages of the strong unstable manifolds.
The leaf $W_f^{u}(x_0)$ is immersed into $\mathbb{T}^3$. In $\mathbb{T}^3$ curves $U(x_i)$ converge to the curve
$U(y_0)$ (dashed curve in the picture). Hence $U(x_0)$ intersects $W_f^{su}(x_0)$ at $y_1$ with
$d_f^{su}(x_0,y_0)\approx d_f^{su}(y_0,y_1)$.} \label{fig6}
\end{figure}
Now we would like to repeat the procedure. Consider another sequence $\{\tilde x_i\in W_f^{wu}(x_0); i\ge 1\}$,
$\tilde x_i\to y_1$ as $i\to\infty$ and corresponding sequence $\{\tilde z_i\in W_f^{wu}(y_0); i\ge 1\}$. Then
$\tilde z_i\to y_2$ as $i\to\infty$, $y_2\in W_f^{su}(x_0)\cap U(x_0)$. In this way by induction we obtain a
sequence of points $\{y_i\in W_f^{su}(x_0)\cap U(x_0); i\ge 1\}$. These points are ordered on $W_f^{su}(x_0)$
--- for any positive $i$ point $y_{i-1}$ lies between $x_0$ and $y_i$. By Claim 2 we know that there are constants
$c_1$ and $c_2$ which depend only on the initial choice of $x_0$ and $y_0$ such that $\forall i\ge 0\;\;\;
c_1>d_f^{su}(y_i,y_{i+1})>c_2$. This guarantees that the set $\{f^{-n}(y_i); n\ge 0, i\ge 0 \}\subset
W_f^{su}(x_0)\cap U(x_0)$ is dense and hence applying Claim 3 one more time we get that $W_f^{su}=U$.
\end{proof}

\subsection{Final remarks}

We did not discuss the proofs of Lemmas 3 and 7. They can be carried out in the same way as the proof of Lemma 5. The technical difficulty with constructing special measure is not present. One can use SRB measures instead (as a matter of fact the construction in Step 2 applied to $W_f^s$ and $W_f^{su}$ will produce SRB measures).

Notice that we used the assumption that $f, g \in \mathcal U$ only to prove Lemmas 1 and 2. So for Theorem 2 we only need to reprove these two lemmas in the new setting. We use a result from~\cite{BI} that states the following.

\begin{BI}
Let $f$ be a partially hyperbolic diffeomorphism of $\mathbb T^3$. Then the lifts of stable and unstable foliations are quasi-isometric and the hence the central distribution is uniquely integrable.
\end{BI}
Thus Lemma 1 is automatic. Proof of Lemma 2 go through with minor differences since we know that $W_f^{su}$ is quasi-isometric.

The bootstrap of regularity of $h$ to the regularity of $f$ and $g$ cannot be done straightforwardly. The reason is the lack of smoothness of weak unstable foliation. Let $N=[\log \lambda_3/\log\lambda_2]$. It is known~\cite{LW} that given $f$ sufficiently $C^1$-close to $L$ the individual leaves of weak unstable foliation are $C^N$ immersed curves. In general the the leaves of weak untable foliation cannot be more than $C^N$ smooth. An example was constructed in~\cite{JPL}. Hence our method cannot lead to smoothness higher than $C^N$.

\end{document}